\documentclass[letterpaper, fancyheadings, 10pt, final]{article}

\def\midformat{
\setlength{\itemsep}{0pt} \setlength{\parindent}{0mm}
\setlength{\parskip}{0.09in} \setlength{\textheight}{180mm}
\setlength{\textwidth}{126mm} \setlength{\evensidemargin}{0in}
\setlength{\oddsidemargin}{0in} \setlength{\topmargin}{0in}
\setlength{\hoffset}{1.0cm} \setlength{\voffset}{2.0cm}
\setlength{\headheight}{0in} \setlength{\headsep}{0in} }
\midformat

\usepackage{psfig}
\usepackage{amssymb,amsmath,amsthm,amsfonts}

\newsavebox{\savepar}

\newtheoremstyle{theorem}% name
 {}% Space above
 {}% Space below
 {\itshape}% Body font
 {}% Indent amount (empty = no indent, \parindent = para indent)
 {\ttfamily}% Thm head font
 {.}% Punctuation after thm head
 {.5em}% Space after thm head: " " = normal interword space;
 {}% Thm head spec (can be left empty, meaning `normal')

\newtheoremstyle{plaintext}% name
 {}% Space above
 {}% Space below
 {\upshape}% Body font
 {}% Indent amount (empty = no indent, \parindent = para indent)
 {\ttfamily}% Thm head font
 {.}% Punctuation after thm head
 {.5em}% Space after thm head: " " = normal interword space;
 {}% Thm head spec (can be left empty, meaning `normal')

\theoremstyle{theorem}
\newtheorem{theorem}{Theorem}[section]
\newtheorem{lemma}[theorem]{Lemma}

\theoremstyle{plaintext}

\begin{document}

\title{\large\bf A new proof of a theorem of Ivanov and Schupp}
\author{Delaram Kahrobaei}
\date{}
\maketitle

\begin{abstract}
The object of this note is to give a very short proof of the
following theorem of Ivanov and Schupp. Let $H$ be a finitely
generated subgroup of a free group $F$ and the index $[F:H]$
infinite. Then there exists a nontrivial normal subgroup $N$ of
$F$ such that $N \cap H = \{1 \}$.
\end{abstract}

As noted in the abstract above the object of this note is to give
a new proof of the following theorem of S.V.Ivanov and P.E.Schupp
[2]
%\cite{IS98}.

\begin{theorem}
Let $H$ be a finitely generated subgroup of a free group $F$ and
the index $[F:H]$ infinite. Then there exists a nontrivial normal
subgroup $N$ of $F$ such that $N \cap H = \{1 \}$.
\end{theorem}

Our proof depends on two facts. The first is a theorem of M. Hall
[1]
%\cite{MH76}
which states that a finitely generated subgroup of a finitely
generated free group is a free factor of a subgroup of finite
index. The second key fact, is the following simple lemma.
\begin{lemma}
\label{lemma} Let $I$ be a normal subgroup of finite index of a
group $F$, and let $L$ be a normal subgroup of $I$. Furthermore
let $b_1 ,..., b_m$ be a complete set of representatives of the
right cosets $Ig$ of $I$ in $F$. Then $$N = \bigcap^m_{i=1}
{b_i}^{-1} L b_i \trianglelefteq F$$
\end{lemma}

\begin{proof} Notice that conjugation induces an automorphism of
$I$ and hence $L^f \trianglelefteq I$ for every $f \in F$. Hence
each ${b_i}^{-1} L b_i \trianglelefteq I$ and therefore $N
\trianglelefteq I$. In fact $N \trianglelefteq F$. To see this,
let $f \in F$. Then $f$ can be written in the form $$f = x b
\text{ } (x \in I, b \in \{b_1,...,b_m\} )$$ Consequently
\begin{enumerate}
\item[(1)] $f^{-1} N f = b^{-1} x^{-1} N x b = b^{-1} N b = b^{-1}
(\bigcap^m_{i=1} {b_i}^{-1} L b_i) b = \bigcap^m_{i=1} (b_i
b)^{-1} L (b_i b)$
\end{enumerate}
Notice that if $b_{i'}$ is the representative of the coset $I b_i
b$ then the map $ i \rightarrow i'$ is a permutation of $\{1,...,
m\}$. Hence we now write $$b_i b = a_i b_{i'} \text{   } (i
=1,...,m, \text{  } a_i \in I)$$ Then it follows from (1) that:
$$f^{-1} N f = \bigcap^m_{i=1} {b_{i'}}^{-1} {a_i}^{-1} L a_i b_{i'}
= \bigcap^m_{i=1} {b_{i'}}^{-1} L b_{i'} = \bigcap^m_{i=1}
{b_i}^{-1} L b_{i}= N$$
\end{proof}

\begin{proof}{\em of the theorem.} If $F$ is not finitely
generated or $H$ is trivial, then the theorem is trivial. So we
assume that $F$ is finitely generated and that $H$ is not trivial.
By the theorem of M.Hall cited above, there exists a subgroup $K
\leq F$ of finite index in $F$ and such that $K$ is the free
product of $H$ a second subgroup $Q$: $$K = H * Q$$ Since $K$ is
of finite index in $F$, the number of conjugates of $K$ in $F$ is
finite and their intersection $I$ is a normal subgroup of $F$ of
finite index (see e.g. A.G. Kurosh [3]%\cite{AGK60}
vol. 1, p. 83)
This can be proved in the same manner that we proved the lemma
above. It follows directly from the proof of the subgroup theorem
for free products due to A.G.Kurosh([3], %\cite{K},
vol. 2, p. 17) that $I \cap H$ is a free factor of $I$:
$$I =(I \cap H) * J$$
Now the normal closure of $J$ in $I$, $L =gp_{I}(J)
\trianglelefteq I$. Observe that by $L \cap (I \cap H) =\{1\}$,
therefore $$L \cap H =\{1\}$$ Now by the lemma $N = \bigcap^m
_{j=1} L^{b_j} \trianglelefteq F$. Note that $N$ is a nontrivial
normal subgroup of $F$, since it contains $[L^{b_1}, L^{b_2},...,
L^{b_m}]$, the subgroup generated by all $m$-fold commutators with
arguments coming out of $L^{b_1}, ..., L^{b_m}$. Then $N$ is the
desired normal subgroup, since $ N \cap H =\{1\}$, using the fact
that $N \leq L$.
\end{proof}

\end{document}